\documentclass{article}

\usepackage{amssymb}
\usepackage{latexsym}
\usepackage{theorem}
\usepackage[all]{xy}

\newenvironment{proof}{\par \noindent{\bf Proof: }}{\hspace{\stretch{1}} $\Box$ \par \mbox{}}
\newcommand{\noproof}{\hspace{\stretch{1}} $\Box$}

\newtheorem{theorem}{Theorem}[section]
\newtheorem{proposition}[theorem]{Proposition}
\newtheorem{lemma}[theorem]{Lemma}
\newtheorem{corollary}[theorem]{Corollary}
{\theorembodyfont{\rmfamily}
\newtheorem{definition}[theorem]{Definition}
\newtheorem{example}[theorem]{Example}
}
\newenvironment{theorem*}{\par \medskip \noindent{\bf Theorem }}{\par \mbox{}}
\newenvironment{lemma*}{\par \medskip \noindent{\bf Theorem }}{\par \mbox{}}

\newcommand{\Hom}{\mathop{\rm Hom}}
\newcommand{\Map}{\mathop{\rm Map}}
\newcommand{\Ob}{\mathop{\rm Ob}}

\newcommand{\supp}{\mathop{\rm Supp}}

\newcommand{\R}{\mathbb{R}}
\newcommand{\C}{\mathbb{C}}
\newcommand{\K}{\mathbb{K}}
\newcommand{\Z}{\mathbb{Z}}
\newcommand{\N}{\mathbb{N}}

\newcommand{\im}{\textrm{im }}

\title{The General Notion of Descent in Coarse Geometry}

\author{Paul D. Mitchener}

\begin{document}

\maketitle

\section*{Abstract}
In this article, we introduce the notion of a functor on coarse spaces being {\em coarsely excisive}- a coarse analogue of the notion of a functor on topological spaces being excisive.  Further, taking cones, a coarsely excisive functor yields a topologically excisive functor, and for coarse topological spaces there is an associated coarse assembly map from the topologically exicisive functor to the coarsely excisive functor.

We conjecture that this coarse assembly map is an isomorphism for uniformly contractible spaces with bounded geometry, and show that the coarse isomorphism conjecture, along with some mild technical conditions, implies that a correspoding equivariant assembly map is injective.  Particular instances of this equivariant assembly map are the maps in the Farrell-Jones conjecture, and in the Baum-Connes conjecture.

\section{Introduction}

In coarse geometry the term {\em descent} usually refers to the fact
that the coarse Baum-Connes conjecture for a group $G$ equipped with
the word length metric, along with a few fairly mild conditions on
the space $EG$, implies rational injectivity of the analytic Novikov
assembly map described in \cite{Kas3,Kas4}.  We refer the reader to
\cite{HR2,HR1,Roe1} for details of the argument.  This descent
principle has been generalised to non-metric coarse structures, such
as those arising from compactifications, as described in \cite{HPR}.

There are similar coarse notions for assembly maps in algebraic $K$-
and $L$-theory arising from controlled topology.  For example, the
descent result involving coarse structures arising from
compactifications was inspired by the work in controlled algebra in
\cite{CP}.  A very strong parallel with methods in the coarse
Baum-Connes conjecture can be found in \cite{Bart}, where methods
used to prove the coarse Baum-Connes conjecture for spaces of finite
asymptotic dimension in \cite{Wright,Yu} are used to show that the
algebraic $K$-theory assembly map (see \cite{BHM,Lod}) is injective
for groups with finite asymptotic dimension.

It is our purpose in this article to give a general theory of coarse
assembly maps and descent techniques that incorporates all of the
above.  We define a general notion of a coarse assembly map.  Our main
result essentially shows that the coarse assembly map being an isomorphism, along with
certain mild technical conditions, implies that the corresponding equivariant coarse assembly map is
injective.

Cases of the equivariant coarse assembly map include the analytic Novikov assembly
map, the Baum-Connes assembly map \cite{BCH}, and the Farrell-Jones assembly map
in algebraic $K$-theory \cite{FJ}.

In principle, geometric arguments used to prove the coarse
Baum-Connes conjecture can be adapted to the new abstract framework.
For example, the results of \cite{Mitch4} immediately tell us
that every coarse assembly map is an isomorphism for
finite coarse cellular complexes.

\section{Coarse Spaces}

Recall that a {\em coarse space} is a set $X$ equipped with a number
of distinguished subsets of the product $X\times X$ called {\em
controlled sets}.

The collection of controlled sets is required to be closed under finite
unions, taking subsets, reflections in the diagonal of $X\times X$,
and composition in the sense that we define
\[ M_1 M_2 = \{ (x,z)\in X\times X \ |\ (x,y)\in M_1 \ (y,z)\in M_2 \textrm{ for some }y\in X
\} \]
for controlled sets $M_1$ and $M_2$.

We further require the diagonal, $\Delta_X = \{ (x,x) \ |\ x\in X
\}$, to be controlled,\footnote{The assumption that the diagonal is
a controlled set is sometimes dropped; a coarse structure where this
axiom does not hold is termed {\em non-unital}.
Many of our definitions should be modified somewhat in the non-unital case; see \cite{Luu}.}
and the union of all controlled sets to be the entire space $X\times
X$.  We refer the reader to \cite{Roe6}, for example, for further
details.

Given a controlled set $M\subseteq X\times X$, and a subset
$S\subseteq X$, we write
$$M[S] = \{ y\in X \ |\ (x,y)\in M \textrm{ for some }x\in S \}$$

For a point $x\in X$, we write $M(x) = M[\{ x \}]$.

If $X$ is a coarse space, and $f,g\colon S\rightarrow X$ are maps
into $X$, the maps $f$ and $g$ are termed {\em close} or {\em
coarsely equivalent} if the set $\{ (f(s) , g(s)) \ |\ s\in S \}$ is
controlled.  We call a subset $B\subseteq X$ {\em bounded} if the
inclusion $B\hookrightarrow X$ is close to a constant map, or
equivalently $B=M(x)$ for some controlled set $M$ and some point
$x\in X$.

The following definition comes from \cite{HPR}.

\begin{definition}
Let $X$ be a Hausdorff space.  A coarse structure on $X$ is said to
be {\em compatible with the topology} if every controlled set is
contained in an open controlled set, and the closure of any bounded
set is compact.

We call a Hausdorff space equipped with a coarse structure that is
compatible with the topology a {\em coarse topological space}.
\end{definition}

Note that any coarse topological space is locally compact.  In such
a space, the bounded sets are {\em precisely} those which are
precompact.  It also follows from the definition that any precompact
subset of the product of the space with itself is an controlled set,
and the closure of any controlled set is an controlled set.

\begin{example}
If $X$ is a proper metric space, the {\em bounded coarse structure}
is the unital coarse structure formed by defining the controlled
sets to be subsets of {\em neighbourhoods of the diagonal}:
$$N_R = \{ (x,y) \in X\times X \ |\ d(x,y)< R \}$$

The bounded sets are simply those which are bounded with respect to
the metric.  If we give the metric space $X$ the bounded coarse
structure, it is a coarse topological space.
\end{example}

The following example is a generalisation of the continuously
controled coarse structure arising from a compactification, as found
for instance in \cite{HPR}.

\begin{example} \label{compactification}
Let $X$ be a coarse topological space, and suppose that $X$ is
contained in a Hausdorff topological space $\overline{X}$ as a
topologically dense subset.  Call the coarse structure already
defined on the space $X$ the {\em ambient coarse structure}.

Write $\partial X = \overline{X}\backslash X$.  Call an open subset
$M\subseteq X\times X$ {\em strongly controlled} if:

\begin{itemize}

\item The set $M$ is controlled with respect to the ambient coarse
structure on $X$.

\item Let $\overline{M}$ be the closure of the set $M$ in the space
$\overline{X}$.  Then
\[ \overline{M} \cap (\overline{X}\times \partial X \cup \partial X
\times \overline{X} ) \subseteq \Delta_{\partial X} .\]

\end{itemize}

Then we define the {\em continuously controlled coarse structure} with respect to
$\overline{X}$ by saying that the controlled sets are composites of subsets
of strongly controlled open sets.

\end{example}

We write $X$ to denote the space $X$ with its ambient coarse
structure, and $X^\mathrm{cc}$ to denote the space $X$ with the new
continuously controlled coarse structure.

\begin{proposition}
The space $X^\mathrm{cc}$ is a coarse topological space.  Further,
if $B\subseteq X^{cc}$ is bounded, then $\overline{B}\cap \partial X
= \emptyset$.
\end{proposition}

\begin{proof}
By construction, we have a coarse structure where every controlled
set is a subset of an open controlled set.

Let $B\subseteq X^\mathrm{cc}$ be bounded.  Then the set $B$ is also
bounded with respect to the ambient coarse structure on the space
$X$.  Since the ambient coarse structure is compatible with the
topology, the closure of the set $B$ in the space $X$ is compact.
Thus the continuously controlled coarse structure is compatible with
the topology.

Since the closure of the set $B$ is compact in the space $X$, it is
also closed in the space $\overline{X}$.  We see
$\overline{B}\subseteq X$, so $\overline{B}\cap \partial X
=\emptyset$.
\end{proof}

Let $X$ and $Y$ be coarse spaces.  Then a map $f\colon X\rightarrow
Y$ is said to be {\em controlled} if for every controlled set
$M\subseteq X\times X$, the image
$$f[M] = \{ (f(x),f(y))\ |\ (x,y)\in M \}$$
is an controlled set.  A controlled map is called {\em coarse} if
the inverse image of any bounded set is bounded.

We can form the category of all coarse spaces and coarse maps.  We
call this category the {\em coarse category}.  We call a coarse map
$f\colon X\rightarrow Y$ a {\em coarse equivalence} if there is a
coarse map $g\colon Y\rightarrow X$ such that the composites $g\circ
f$ and $f\circ g$ are close to the identities $1_X$ and $1_Y$
respectively.

Coarse spaces $X$ and $Y$ are said to be {\em coarsely equivalent},
and we write $X\sim Y$, if there is a coarse equivalence between
them.

\begin{definition}
Let $X$ and $Y$ be coarse spaces.  Then we define the product,
$X\times Y$ to be the Cartesian product of the sets $X$ and $Y$
equipped with the coarse structure defined by saying a subset
$M\subseteq (X\times Y)\times (X\times Y)$ is controlled if it is
the subset of a set of the form
\[ \{ (u,v,x,y) \ |\ (u,x)\in M_1,\ (v,y)\in M_2 \} \]
where $M_1\subseteq X\times X$ and $M_2\subseteq Y\times Y$ are
controlled sets.
\end{definition}

Note that the above product is not in general a product in the
category-theoretic sense, since the projections onto the factors are
not coarse maps.

\begin{definition}
Let $X$ be a coarse space, let $\sim$ be an equivalence relation on
$X$, and let $X/\sim$ be the set of equivalence classes.  Let $\pi
\colon X\rightarrow X/\sim$ be the quotient map sending each point
$x\in X$ to its equivalence class, $\pi (x)$.

We define the {\em quotient coarse structure} on $X/\sim$ by saying
a subset $M\subseteq X/\sim \times X/\sim$ is controlled if and only
if the union $M = \pi [M']$ for some controlled set $M'\subseteq
X\times X$.
\end{definition}

Note that the quotient map $\pi \colon X\rightarrow X/\sim$ is not
in general a coarse map.

\begin{definition}
Let $\{ X_i \ |\ i\in I \}$ be a collection of coarse spaces.  Then,
as a set, the {\em coarse disjoint union}, $\vee_{i\in I} X_i$ is
the disjoint union of the sets $X_i$.

A subset $M\subseteq (\vee_{i\in I} X_i)\times (\vee_{i\in
I} X_i)$ controlled if it is a subset of a union of the form
\[ \left( \bigcup_{i\in I} M_i \right) \cup
\left( \bigcup_{i,j\in I} B_i \times B_j \right) \]
where each set $M_i\subseteq X_i \times X_i$ is controlled, and
$B_j\subseteq X_j$ is bounded.
\end{definition}

For example, let $\R_+$ be the space $[0,\infty )$ equipped with the
bounded coarse structure arising from the usual metric.  Then $\R_+
\vee \R_+ = \R$, where the real line, $\R$, is again given the
bounded coarse structure.

\section{Homotopy}

The following definition comes from \cite{Mitch4}.

\begin{definition}
Let $R$ be a the topological space $[0,\infty )$ equipped with a coarse structure
compatible with the topology.  We call the space $R$ a {\em
generalised ray} if the following conditions hold.

\begin{itemize}

\item Let $M,N\subseteq R\times R$ be controlled sets.  Then the sum
\[ M+N = \{ (u+x,v+y)\ |\ (u,v)\in M,\ (x,y)\in N \} \]
is controlled.

\item Let $M\subseteq R\times R$ be a controlled set.  Then the set
\[ M^s = \{ (u,v)\in R\times R\ |\ x\leq u,v\leq y,\ (x,y)\in M \}
\]
is controlled.

\item Let $M\subseteq R\times R$ be a controlled set, and $a\in R$.
Then the set
\[ a+M = \{ (a+x,a+y)\ |\ (x,y)\in M \} \]
is controlled.

\end{itemize}

\end{definition}

For example, the space $\R_+$ (with the bounded coarse structure) is a
generalised ray.  The space $[0,\infty )$ equipped with the
continuously controlled coarse structure arising from the one point
compactification is also a generalised ray.

In order to look at the notion of homotopy for coarse spaces,
we first consider cylinders. Our first definition is inspired by section 3 of \cite{Dran}.

\begin{definition}
Let $X$ be a coarse space equipped with a coarse map
$p\colon X\rightarrow R$.  Then we define the {\em $p$-cylinder of $X$}:
$$I_p X = \{ (x,t)\in X\times R \ |\ t\leq p(x)+1 \}$$

We define coarse maps $i_0,i_1\colon X\rightarrow I_pX$ by the
formulae $i_0 (x) = (x,0)$ and $i_1 (x)=(x,p(x)+1)$ respectively.
\end{definition}

The following definition is inspired by \cite{MiS}.

\begin{definition}
Let $f_0,f_1\colon X\rightarrow Y$ be coarse maps.  An {\em elementary coarse homotopy}
between $f_0$ and $f_1$ is a coarse map
$H\colon I_p X\rightarrow Y$ for some $p\colon X\rightarrow R$ such that $f_0 = H\circ i_0$
and $f_1 = H\circ i_1$.

More generally, we call the maps $f_0$ and $f_1$ {\em coarsely
homotopic} if they can be linked by a chain of elementary coarse
homotopies.
\end{definition}

A coarse map $f\colon X\rightarrow Y$ is termed a {\em coarse
homotopy equivalence} if there is a coarse map $g\colon Y\rightarrow
X$ such that the compositions $g\circ f$ and $f\circ g$ are coarsely
homotopic to the identities $1_X$ and $1_Y$ respectively.

\begin{example}
Let $X$ and $Y$ be spaces, and let $p\colon X\rightarrow R$ be a
coarse map.  Consider two close coarse maps $f_0\colon X\rightarrow
Y$ and $f_1\colon X\rightarrow Y$.  Then we
can define a coarse homotopy $H\colon I_pX\rightarrow Y$ between the
maps $f_0$ and $f_1$ by the formula
$$H(x,t) = \left\{ \begin{array}{ll}
f_0 (x) & x< 1 \\
f_1 (x) & x\geq 1 \\
\end{array} \right.$$
\end{example}

Thus, close maps are also coarsely homotopic.  In particular, any
coarse equivalence is also a coarse homotopy equivalence.

\section{Coarse Homology}

The following definition comes from \cite{HRY}, where it is a condition for the existence of a
coarse version of the Mayer-Vietoris sequence.

\begin{definition}
Let $X$ be a coarse space.  Then we call a decomposition $X=A\cup B$
{\em coarsely excisive} if for every controlled set $m\subseteq X\times X$ there is a controlled set
$M\subseteq X\times X$ such that $m(A)\cap m(B)\subseteq M(A\cap B)$.
\end{definition}

The following definition is now a slight variant of the definition in \cite{Mitch4}.

\begin{definition}
A {\em coarse homology theory} is a collection of functors,
$\{ h^\mathrm{coarse}_n\ |\ n\in \Z \}$, from the coarse category to the category of abelian groups such
that the following axioms are satisfied:

\begin{itemize}

\item Let $X$ be a coarse space,
and let $f,g\colon X\rightarrow Y$ be coarsely homotopic coarse
maps.  Then the induced homomorphisms $f_\ast , g_\ast \colon
h^\mathrm{coarse}_\ast (X)\rightarrow h^\mathrm{coarse}_\ast (Y)$
are equal.

\item Let $X$ be a coarse space, and let $X=A\cup B$ be a coarsely excisive decomposition.  Let $i\colon A\cap B\hookrightarrow A$, $j\colon A\cap B\hookrightarrow B$, $k\colon A\hookrightarrow X$, and $l\colon B\hookrightarrow X$ be the associated inclusion maps.  Then there are natural maps $\partial \colon h^\mathrm{coarse}_n (X) \rightarrow h^\mathrm{coarse}_{n-1} (X)$ such that we have a long exact sequence
$$\rightarrow h^\mathrm{coarse}_n(A\cap B)\stackrel{(i_\ast , -j_\ast )}{\rightarrow}h^\mathrm{coarse}_n (A)\oplus h^\mathrm{coarse}_n (B) \stackrel{k_\ast + l_\ast}{\rightarrow} h^\mathrm{coarse}_n (X) \stackrel{\partial}{\rightarrow} h^\mathrm{coarse}_{n-1}(A\cap B)\rightarrow$$

\end{itemize}

The above long exact sequence is called a {\em coarse Mayer-Vietoris sequence}.

A coarse homology theory $\{ h^\mathrm{coarse}_n \ |\ n\in \Z \}$ is said to satisfy the
{\em disjoint union axiom} if the following condition also holds.

\begin{itemize}

\item Let $\{ X_i \ |\ i\in I \}$ be a family of coarse spaces, and let $j_i \colon X_i \rightarrow \vee_{i\in I} X_i$ be the canonical inclusion.  Then the map
$$\oplus_{i\in I} (j_i)_\ast \colon \oplus_{i\in I} h^\mathrm{coarse}_n (X_i)\rightarrow h^\mathrm{coarse}_n \left( \vee_{i\in I} X_i \right)$$
is an isomorphism for all $n\in \Z$.

\end{itemize}

\end{definition}

The above axioms are a coarse variant of the
Eilenberg-Steenrod axioms used to define a generalised homology
theory in the world of topological spaces.  There is a corresponding
disjoint union axiom in the topological world.

The following definition essentially comes from \cite{HPR,Roe1}.

\begin{definition}
We call a coarse space $X$ {\em flasque} if there is a map
$t\colon X\rightarrow X$ such that:

\begin{itemize}

\item Let $B\subseteq X$ be bounded.  Then there exists $N\in \N$
such that $t^n[X]\cap B = \emptyset$ for all $n\geq N$.

\item Let $M\subseteq X\times X$ be controlled.  Then the union
$\bigcup_{n\in \N} t^n [M]$ is controlled.

\item The map $t$ is close to the identity map.

\end{itemize}

\end{definition}

In the above definition, the map $t$ is controlled by the second property, but is not necessarily coarse.

Observe that any generalised ray is flasque; the relevant map $t\colon
R\rightarrow R$ is defined by the formula $t(s) =s+1$.

\begin{definition}
We call a functor, $E$, from the coarse category to the category of
spectra {\em coarsely excisive} if the following conditions hold.

\begin{itemize}

\item The spectrum $E(X)$ is weakly contractible whenever the coarse
space $X$ is flasque.

\item The functor $E$ takes coarse homotopy equivalences to weak homotopy
equivalences of spectra.

\item For a coarsely excisive decomposition $X=A\cup B$ we have a homotopy
push-out diagram
\[ \begin{array}{ccc}
E(A\cap B) & \rightarrow & E(A) \\
\downarrow & & \downarrow \\
E(B) & \rightarrow & E(X) \\
\end{array} . \]

\item Up to homotopy, the functor $E$ takes disjoint unions in the coarse category to coproducts
in the category of spectra.

\end{itemize}

\end{definition}

The following is immediate.

\begin{proposition}
Let $E$ be a coarsely excisive functor.  Then the sequence of
functors $X\mapsto \pi_n E(X)$ forms a coarse homology theory
satisfying the disjoint union axiom.
\noproof
\end{proposition}

\begin{definition}
Let $X$ be a compact Hausdorff space.  Then we define the {\em open
cone} to be the space
\[ {\mathcal O}X = \frac{X\times [0,1)}{X\times \{ 0 \} } \]
equipped with the continuously controlled coarse structure arising
from the compactification
\[ \frac{X \times [0,1]}{X\times \{ 0 \} } \]
\end{definition}

Suppose that $X$ is a subset of the unit sphere of some Hilbert space, $H$.
Let $\varphi \colon [0,1)\rightarrow [0,\infty )$ be a homeomorphism.
Then we have an induced map $\varphi_\ast \colon {\mathcal O}X\rightarrow H$ defined by writing
$$\varphi_\ast (x,t) = \varphi (t)x$$
and we define ${\mathcal O}_\varphi X$ to be the image of the above
map $\varphi_\ast$, equipped with the bounded coarse structure
arising from the metric of the Hilbert space $H$.

Proposition 6.2.1 of \cite{HR1} tells us the following.

\begin{proposition} \label{conrel}
Let $X$ be a compact subset of the unit sphere of a Hilbert space, $H$.  Then:

\begin{itemize}

\item Let $\varphi \colon [0,1)\rightarrow [0,\infty )$ be a homeomorphism.
Then every controlled set for the space ${\mathcal O}_\varphi X$ is also controlled for the space
${\mathcal O} X$.

\item Let $M$ be a controlled set for the space ${\mathcal O}X$.
Then there is a homeomorphism $\varphi \colon [0,1)\rightarrow [0, \infty )$
such that the set $M$ is controlled for the space ${\mathcal O}_\varphi X$.

\end{itemize}

\noproof
\end{proposition}

The idea of a functor being coarsely excisive is an analogue of the
following notion from \cite{WW}.

\begin{definition}
Let $\mathcal C$ be the category of topological spaces that are
homotopy-equivalent to $CW$-complexes.  We call a functor, $F$, from
the category $\mathcal C$ to the category of spectra {\em excisive} if it takes homotopy equivalences of spaces to weak homotopy
equivalences of spectra, up to homotopy takes disjoint unions of topological spaces
to coproducts in the category of spectra, and given a space
$X=U\cup V$, where $U$ and $V$ are open sets we have a homotopy
push-out diagram
\[ \begin{array}{ccc}
F(U\cap V) & \rightarrow & F(U) \\
\downarrow & & \downarrow \\
F(V) & \rightarrow & F(X) \\
\end{array} . \]
\end{definition}

Given a strongly excisive functor $F$, the sequence of functors $X\mapsto \pi_n F(X)$ forms a generalized
homology theory satisfying the disjoint union axiom.

\begin{theorem} \label{CH}
Let $E$ be a coarsely excisive functor.  Then the assignment $X\mapsto E({\mathcal O}X)$ defines a
strongly excisive functor on the category of compact
Hausdorff spaces homotopy-equivalent to $CW$-complexes.
\end{theorem}

\begin{proof}
Let $f\colon X\rightarrow Y$ be a continuous map.
Then the obvious induced map $f_\ast \colon {\mathcal O}X \rightarrow {\mathcal O}Y$ is coarse.

The space $[0,1)$, equipped with the continuously controlled coarse structure
arising from the one point compactification, can be considered to be
a generalised ray.  We have a coarse map $p\colon {\mathcal
O}X\rightarrow [0,1]$ defined by the formula $p[(x,t)] =t$.

Consider a homotopy $F\colon X\times [0,1]\rightarrow Y$.
Then the induced map is a coarse homotopy $F\colon I_p {\mathcal O}X \rightarrow
Y$.  Thus, given a homotopy equivalence $f\colon X\rightarrow Y$, we
have an induced homotopy equivalence $f_\ast \colon E({\mathcal
O}X)\rightarrow E({\mathcal O}Y)$.

Suppose we have a decomposition $X=U\cup V$, where $U$ and $V$ are open subsets.
Embed $X$ in the unit sphere of a Hilbert space, $H$, and choose a homeomorphism
$\varphi \colon [0,1)\rightarrow [0,\infty)$.
Since the space $X$ is precompact, and the subsets $U$ and $V$ are open, for any $R>0$ we can find $S>0$
such that $N_R [{\mathcal O}_\varphi U]\cap N_R [{\mathcal
O}_\varphi V] \subseteq N_S [{\mathcal O}_\varphi (U\cap V)]$.

It follows that the decomposition
${\mathcal O}_\varphi X = {\mathcal O}_\varphi U \cup {\mathcal O}_\varphi V$
is coarsely excisive.  It follows by proposition \ref{conrel} that
the decomposition ${\mathcal O}X = {\mathcal O}U \cup {\mathcal O}V$
is also coarsely excisive.  It follows that we have a homotopy
push-out square
\[ \begin{array}{ccc}
E({\mathcal O}(U\cap V)) & \rightarrow & E({\mathcal O}U) \\
\downarrow & & \downarrow \\
E({\mathcal O}V) & \rightarrow & E({\mathcal O}X) \\
\end{array} . \]

Let $\{ X_i \ |\ i\in I \}$ be a collection of locally compact Hausdorff spaces.  Then
$${\mathcal O}(\amalg_{i\in I} X_i ) = \vee_{i\in I} {\mathcal O} X_i$$

Thus the functor $E{\mathcal O}$ takes disjoint unions of topological
spaces to colimits of spectra, and we are done.
\end{proof}

\section{Coarse Assembly}

Let $X$ be a coarse topological space.  Then we define the {\em open square} to be the space
${\mathcal S}X = X\times [0,1)$ equipped with the continuously
controlled coarse structure arising by considering ${\mathcal S}X$
as a dense subspace of the topological space $X\times [0,1]$.

We define the {\em closed cone} to be the quotient coarse space
${\mathcal C}X = X\times [0,1]/X\times \{ 0 \}$.  Note that there is
no element of continuous control in the definition of the coarse
structure on the closed cone; the coarse structure comes from the
ambient coarse structure on the space $X$, and taking products and
quotients.

\begin{proposition} \label{sxf}
Let $X$ be a coarse topological space.  Then the open square
${\mathcal S}X$ is flasque.
\end{proposition}

\begin{proof}
Define a map $\alpha \colon [0,1)\rightarrow [0,1)$ by the formula
\[ \alpha (s) = \frac{1}{2-s} . \]

Then the map $\alpha$ is a continuous monotonic increasing map, and
for any point $s_0 \in [0,1)$, if we define the sequence $(s_n)$
iteratively by the formula $s_{n+1} = \alpha (s_n)$, then
$\lim_{n\rightarrow \infty} s_n = 1$.

Define a map $t\colon {\mathcal S}X \rightarrow {\mathcal S}X$ by
the formula
\[ t(x,s) = (x, \alpha (s)) . \]

Now, let $M\subseteq {\mathcal S}X \times {\mathcal S}X$ be strongly
controlled.  Then, looking at the ambient coarse structure, we have
\[ M \subseteq M_X \times [0,1) \times [0,1) \]
where $M_X\subseteq X\times X$ is controlled.  Here we are indulging
in some slight abuse of notation involving the order of our factors.

Further,
\[ \overline{M}\cap ((X\times [0,1] \times X\times \{ 1 \} ) \cup ( X\times \{ 1
\} \times (X\times [0,1])) \subseteq X_{X\times \{ 1 \} } . \]

It follows that
\[ t[M] \subseteq M_X \times [0,1)\times [0,1) \]
and
\[ t[\overline{M}] \cap ((X\times [0,1] \times X\times \{ 1 \} ) \cup ( X\times \{ 1
\} \times (X\times [0,1])) \subseteq X_{X\times \{ 1 \} } . \]

Iterating, we see that $\cup_{n\in \N }t^n [M]$ is controlled with
respect to the continuously controlled coarse structure on the space
${\mathcal S}X$.  In particular, it follows that the map $t$ is
coarse.

Let
\[ M = \{ ((x,s),(x,\alpha (s)) \ |\ x\in X,\ s\in [0,1) \} . \]

Then the set $M$ is certainly controlled with respect to the ambient
coarse structure on $X$.

Let $(s_n)$ be a sequence in the space $[0,1)$.  Observe that $s_n
\rightarrow 1$ as $n\rightarrow \infty$ if and only if $\alpha
(s_n)\rightarrow 1$ as $n\rightarrow \infty$.  Hence
\[ \begin{array}{rcl}
 \overline{M} \cap ((X\times [0,1] \times X\times \{ 1 \} ) \cup ( X\times \{ 1
\} \times (X\times [0,1])) & = & \{ ((x,1),(x,1))\ |\ x\in X \} \\
& = &\Delta_{X\times \{ 1\} } \\
\end{array}
. \]

We see that the set $M$ is controlled with respect to the
continuously controlled coarse structure on the space ${\mathcal
S}X$.  Hence the map $t$ is close to the identity $1_{{\mathcal
S}X}$.

Fianlly, let $B\subseteq {\mathcal S}X$ be bounded.  Then, as we
remarked earlier, $B\subseteq X\times [0,a]$ for some $a\in [0,1)$.
It follows that we have $N\in \N$ such that $t^n [{\mathcal S}X]\cap
B = \emptyset$ for all $n\geq N$.

In conclusion, we see that the space ${\mathcal S}X$ is flasque.
\end{proof}

\begin{proposition}
Let $X$ be a flasque coarse space.  Then the closed cone ${\mathcal
C}X$ is also flasque.
\end{proposition}

\begin{proof}
Let $t\colon X\rightarrow X$ be a coarse map with the required
properties to make the space $X$ flasque.  Let $\pi \colon X\times
[0,1]\rightarrow {\mathcal C}X$ be the quotient map.

Define a map $\tilde{t}\colon {\mathcal C}X\rightarrow {\mathcal
C}X$ by the formula
\[ \tilde{t} ([x,s]) = [t(x),s] . \]

Let $\tilde{M}\subseteq {\mathcal C}X\times {\mathcal C}X$ be a
controlled set.  Then $M\subseteq \pi [M\times [0,1]\times [0,1]]$
where $M$ is a controlled set for the coarse space $X$; here we are
indulging in some mild abuse of notation involving the order of
factors.  We know that the set $\bigcup_{n\in \N} t^n[M]$ is
controlled.

It follows that the set $\bigcup_{n\in \N} \tilde{t}^n[\tilde{M}]$
is controlled.   In particular, we have shown that the map
$\tilde{t}$ is controlled.

Now, the map $t$ is close to the identity map on the space $X$.
Therefore the set
\[ M = \{ (x,t(x))\ |\ x\in X \} \]
is controlled.  Let
\[ \tilde{M} = \{ ([x,s],\tilde [x,s] ) \ |\ [x,s]\in {\mathcal C}X
\} . \]

Observe
\[ \tilde{M} = \{ ( [x,s],[t(x),s] ) \ |\ [x,s]\in {\mathcal C}X \}
\subseteq \pi [M \times [0,1] ] \]
so the set $\tilde{M}$ is also controlled, and the map $\tilde{t}$
is close to the identity map.

Finally, let $\tilde{B}\subseteq {\mathcal C}X$ be bounded.  As above, we
have $\tilde{B}\subseteq \pi [B\times [0,1]]$ for some bounded set
$B\subseteq X$.  Now we have a natural number $N$ such that
$t^n[X]\cap B = \emptyset$ whenever $n\geq N$.  Hence
$\tilde{t}^n[{\mathcal C}X] \cap \tilde{B} = \emptyset$ whenever
$n\geq N$.

Thus the space ${\mathcal C}X$ is also flasque, as claimed.
\end{proof}

\begin{lemma} \label{csefib}
Let $X$ be a Hausdorff space equipped with a coarse structure
compatible with the topology.  Let $E$ be a coarsely excisive
functor.
Then we have a weak fibration
\[ E(X) \stackrel{j}{\rightarrow} E({\mathcal C}X) \stackrel{v}{\rightarrow} E({\mathcal O}X) \]
\end{lemma}

\begin{proof}
We have a coarsely excisive decomposition
\[ {\mathcal O}X = \frac{X\times [0,1/2]}{X\times \{ 0 \} } \cup X
\times \left[ \frac{1}{2} ,1 \right) . \]

Observe that the first space in the decomposition is coarsely
equivalent to the space ${\mathcal C}X$, and the second space is
coarsely equivalent to the space ${\mathcal S}X$.  By proposition
\ref{sxf}, the space ${\mathcal S}X$ is flasque, so the space
$E({\mathcal S}X)$ is weakly contractible.

Hence, by the homotopy push-out property of the functor $E$, we have a weak fibration
\[ E(X) \rightarrow E({\mathcal C}X) \rightarrow E({\mathcal O}X) \]
as claimed.
\end{proof}

\begin{definition}
We call the boundary map
\[ \partial \colon \Omega E({\mathcal O}X ) \rightarrow E(X) \]
associated to the above weak fibration the {\em coarse assembly map}
associated to the functor $E$.
\end{definition}

\begin{definition}
Let $X$ be a coarse space.  We say $X$ has {\em bounded geometry} if
it is coarsely equivalent to a space $Y$ such that for any
controlled set $M$, there is a number $k$ such that $|M(x)|\leq k$
for all $x\in Y$.
\end{definition}

For example, any subset of Euclidean space, $\R^n$, has bounded
geometry.

\begin{definition}
Let $X$ be a coarse topological space.  We call $X$ {\em uniformly
contractible} if for every open controlled set $M$, there is a
controlled set $N\supseteq M$ such that for all
$x\in X$, the inclusion $M(x)\hookrightarrow
N(x)$ is homotopic to the constant map onto the point $x$.
\end{definition}

For example (see \cite{Roe6}), a contractible metric space on which some
group acts cocompactly by isometries is uniformly contractible.

\begin{definition}
The {\em coarse isomorphism conjecture} associated to the functor
$E$ asserts that the coarse assembly map is a stable equivalence
whenever $X$ is a uniformly contractible space with bounded
geometry.
\end{definition}

\begin{proposition} \label{wcx}
Let $E$ be a coarsely excisive functor.  The coarse isomorphism
conjecture holds for a space $X$ if and only if the spectrum
$E({\mathcal C}X)$ is weakly contractible.
\end{proposition}

\begin{proof}
We have a long exact sequence
\[ \pi_{n+1}E({\mathcal C}X) \stackrel{v_\ast}{\rightarrow}
\pi_{n+1} E({\mathcal O}X) \stackrel{\partial_\ast}{\rightarrow}
\pi_n E(X) \stackrel{j_\ast}{\rightarrow} \pi_n E({\mathcal C}X ) \]

If the spectrum $E({\mathcal C}X)$ is weakly contractible, then the
groups $\pi_{n+1}E({\mathcal C}X)$ and $\pi_n E({\mathcal C}X)$ are
both zero, so the map $\partial_\ast$ is an isomorphism.  It follows
that the coarse isomorphism conjecture holds for $X$.

Conversely, suppose the coarse isomorphism conjecture holds for $X$.
Then the above map $\partial_\ast$ is an isomorphism, and the maps
$j_\ast$ and $v_\ast$ are both zero.

But we can extend the above long exact sequence to the right:
\[ \pi_n E(X) \stackrel{j_\ast}{\rightarrow} E({\mathcal C}X)
\stackrel{v_\ast}{\rightarrow} \pi_n E({\mathcal O}X) \]

We know that $j_\ast =0$, $v_\ast =0$, and $\im j_\ast = \ker
v_\ast$.  Thus $\pi_n E({\mathcal C}X) =0$, and the spectrum
$E({\mathcal C}X)$ is weakly contractible, as desired.
\end{proof}

\section{Examples}

Let $R$ be a ring.  We call a category, $\mathcal A$, an
{\em $R$-algebroid} (see \cite{Mi}) if each morphism set $\Hom
(a,b)_{\mathcal A}$ is a left $R$-module, and composition of
morphisms
\[ \Hom (b,c)_{\mathcal A}\times \Hom (a,b)_{\mathcal A}\rightarrow
\Hom (a,c)_{\mathcal A} \]
is $R$-bilinear.

Let $\mathcal A$ be an $R$-algebroid, and consider objects
$a,b\in \Ob ({\mathcal A})$.  Then an object $a\oplus b$ is called a
{\em biproduct} of the objects $a$ and $b$ if it comes equipped with
morphisms $i_a \colon a\rightarrow a\oplus b$, $i_b \colon
b\rightarrow a\oplus b$, $p_a \colon A\oplus B\rightarrow A$, and
$p_B \colon a\oplus b \rightarrow b$ satisfying the equations
$$p_a i_a = 1_a \qquad p_b i_b = 1_b \qquad i_a p_a + i_b p_b = 1_{a\oplus b}$$

Observe that a biproduct is simultaneously a product and a coproduct.  An $R$-algebroid $\mathcal R$
is called {\em additive} if it has a zero object (that is, and
object that is simultaneously initial and terminal), and every pair
of objects has a biproduct.

The following construction can be found in several articles on controlled algebraic $K$-theory; see for example \cite{ACFP,Bart,CP}.  Our functorial approach comes from \cite{Weiss}.

\begin{definition}
Let $X$ be a metric space, and let $\mathcal A$ be an additive $R$-algebroid.  Then a {\em geometric $\mathcal A$-module} over $X$ is a functor, $M$, from the collection of bounded subsets of $X$\footnote{Regarded as a category by looking at the usual partial ordering.} to the category $\mathcal A$ such that for any bounded set $B$ the natural map
\[ \oplus_{x\in B} M(\{ x \}) \rightarrow M(B) \]
induced by the various inclusions is an isomorphism, and the {\em support}
$$\supp (M) = \{ x\ |\ M(\{ x\}) \neq 0 \}$$
is a locally finite subset of $X$.
\end{definition}

For convenience, we write $M_x = M(\{ x \})$.  A morphism $\phi \colon M\rightarrow N$ between geometric $\mathcal A$-modules over $X$ is a collection of
morphisms $\phi_{x,y} \colon M_y \rightarrow N_x$ in the category
$\mathcal A$ such that for each fixed point $x\in X$, the morphism $\phi_{x,y}$ is non-zero for only finitely many points
$y\in X$, and for each fixed point $y\in X$, the morphism $\phi_{x,y}$ is non-zero for only finitely many points $x\in X$.

Composition of morphisms $\phi \colon M\rightarrow N$ and $\psi \colon N\rightarrow P$ is defined by the formula
\[ (\psi \circ \phi )_{x,y} (\eta ) = \sum_{z\in X} \psi_{x,z} \circ \phi_{z,y} (\eta ) . \]

The local finiteness condition ensures this makes sense.  We define the {\em support} of a morphism $\phi$
\[ \supp (\phi ) = \{ (x,y)\in X \ |\ \phi_{x,y}\neq 0 \} . \]

\begin{definition}
The category ${\mathcal A}(X)$ consists of all geometric
$\mathcal A$-modules over $X$ and morphisms such that the support is
controlled with respect to the coarse structure of $X$.
\end{definition}

Observe that ${\mathcal A}(X)$ is again an additive $R$-algebroid.
In particular, we can form its algebraic $K$-theory spectrum ${\mathbb K}
{\mathcal A}(X)$ (constructed for instance in \cite{PW}).  The
algebraic $K$-theory groups $K_n {\mathcal A}(X)$ are the stable
homotopy groups of this spectrum.

Given a coarse map $f\colon X\rightarrow Y$, we have an induced additive functor 
$f_\ast \colon {\mathcal A}(X)\rightarrow {\mathcal A}(Y)$.

For a geometric $\mathcal A$-module, $M$, over $X$, we define a geometric $\mathcal A$-module $f_\ast [M]$ by writing $f_\ast [M](S) = M(f^{-1}[S])$.

Let $\phi \colon M\rightarrow N$ be a morphism in the category ${\mathcal A}(X)$.  We define a morphism $f_\ast [\phi ] \colon f_\ast [M]\rightarrow f_\ast [N]$ by writing
\[ f_\ast [\phi ]_{h_1,h_2} = \sum_{g_1\in f^{-1}(h_1) \atop g_2\in f^{-1}(h_2)} \phi_{g_1,g_2} .\]

With these induced maps, we have a functor $X\mapsto {\mathcal A}(X)$, and so a functor  $X\mapsto {\mathbb K}{\mathcal A}(X)$. 

\begin{theorem} \label{Ahom}
The functor $X\mapsto {\mathbb K}{\mathcal A}(X)$ is coarsely excisive.
\end{theorem}

\begin{proof}
To begin with, let $X$ be a coarse space, and let $p\colon
X\rightarrow R$ be a coarse map to a generalised ray.  Consider the inclusion $i_0 \colon X\rightarrow I_pX$
defined by the formula $i(x)=(x,0)$.  Let $q \colon I_p X\rightarrow
X$ be defined by the formula $q(x,t) = x$.
Then certainly $q\circ i_0 = 1_X$, and $i_0\circ q (x,t) = (x,0)$.

Let $M$ be a geometric $\mathcal A$-module over $I_p X$.  Then
$$(i_0 \circ q)_\ast (E)_{(x,t)} =\left\{ \begin{array}{ll}
0 & t\neq 0 \\
\oplus_{0\leq t\leq p(x)+1} E_{(x,t)} & t=0 \\
\end{array} \right.$$

It follows that we have a natural isomorphism $M\rightarrow (i_0\circ q)_\ast (M)$.
Hence the functor $(i_0)_\ast \colon {\mathcal A}[X]\rightarrow {\mathcal A}[I_pX]$ is an equivalence of categories.
The same argument proves that the functor
$(i_1)_\ast \colon {\mathcal A}[X]\rightarrow {\mathcal A}[I_p X]$ where $i_1 (x) = (x,p(x)+1)$ is also an equivalence of categories.

Certainly an equivalence of categories induces a homotopy equivalence at
the level of $K$-theory spectra.  We conclude that the functor $X\mapsto \K {\mathcal A}[X]$ takes
coarse homotopy equivalences to homotopy equivalences.  In
particular, if two maps $f,g\colon X\rightarrow Y$ are coarsely
equivalent, the maps $f_\ast , g_\ast \colon \K {\mathcal
A}[X]\rightarrow \K {\mathcal A}[Y]$ are coarsely homotopic.

Now, let $X$ be flasque.  Then we have a map $\tau \colon X\rightarrow X$
which gives us the flasqueness property.  Given a geometric
$\mathcal A$-module, $M$, the direct sum
\[ \oplus_{n=0}^\infty \tau^n_\ast [M] \]
is also a geometric $\mathcal A$-module since the powers of $\tau$
eventually leave any bounded subset of $X$.

Given a morphism $f\colon M\rightarrow N$ in the category ${\mathcal A}[X]$, the fact that the union
$\bigcup_{n=0}^\infty \tau^n[M]$ is controlled whenever $M\subseteq
X\times X$ is a controlled subset tells us that we have an induced
morphism
\[ \sum_{n=0}^\infty f_\ast \tau_\ast^n \colon \oplus_{n=0}^\infty \tau^n_\ast [M]
\rightarrow \oplus_{n=0}^\infty \tau^n_\ast [N] \]
in the category ${\mathcal A}[X]$.

Hence, by definition of the $K$-theory of an additive category (see
for example \cite{Wa1}), given a $K$-theory class $\lambda \in K_n
{\mathcal A}[X]$, we have an induced $K$-theory class
\[ \sum_{n=0}^\infty f_\ast \tau_\ast^n [\lambda] = \lambda +
\tau_\ast [\lambda ] + \tau^2_\ast [\lambda ] + \cdots \in \K {\mathcal A}[X] . \]

Now the map $\tau$ is close to the identity map on $X$.  Hence, by
the first part of the current proof, the map $\tau_\ast$ is the
identity at the level of $K$-theory groups.  We see that
\[ \sum_{n=0}^\infty f_\ast \tau_\ast^n [\lambda] = \lambda +
\lambda + \cdots = \lambda + (\lambda + \lambda + \cdots ) . \]

Thus, by an Eilenberg swindle, we see that $\lambda =0$.  Thus $K_n
{\mathcal A}[X] =0$, meaning the stable homotopy groups of the
spectrum $\K {\mathcal A}[X]$ are all zero.  Therefore the spectrum
spectrum $\K {\mathcal A}[X]$ is weakly contractible.

We now look at homotopy push-out squares.  let $X = U\cup V$ be a
coarsely excisive decomposition of a coarse space $X$.  Consider the
sequence
$$0\rightarrow {\mathcal A}[U\cap V]\stackrel{(i_\ast, j_\ast)}{\rightarrow} {\mathcal A}[U]\oplus {\mathcal A}[V] \stackrel{k_\ast - l_\ast}{\rightarrow} {\mathcal A}[X]\rightarrow 0$$
where $i,j,k,l$ are the relevant inclusions.

Clearly the functor $(i_\ast ,j_\ast )$ is faithful, the functor $k_\ast - l_\ast$ is surjective on each morphism set, and $(k_\ast -l_\ast) \circ (i_\ast , j_\ast ) =0$.

We claim that the above sequence is a short exact sequence.
Consider morphisms $\phi ,\psi \colon M\rightarrow N$ between
geometric ${\mathcal A}$-modules over $U$ and $V$ respectively such
that $k_\ast \phi - l_\ast \psi =0$, and the supports of $\phi$ and
$\psi$ are controlled.

We see that $\phi_{x,y}=0$ and $\psi_{x,y}=0$ unless $x,y\in U\cap
V$, and $\phi_{x,y}=\psi_{x,y}$ when $x,y\in U\cap V$.  Define a
morphism $\theta \colon M|_{U\cap V}\rightarrow N|_{U\cap V}$ by
restricting the morphism $\phi$.    Then certainly $(\phi, \psi) =
(i_\ast , j_\ast )(\theta )$.

The support of the morphism $\theta$ is certainly controlled.  By the fibration theorem in algebraic $K$-theory (see
\cite{McC,Wa1}). we have a fibration
\[ \K {\mathcal A}[U\cap V] \rightarrow \K {\mathcal A}[U]\vee \K {\mathcal
A}[V] \rightarrow \K {\mathcal A}[X] . \]

It follows that we have a homotopy push-out
\[ \begin{array}{ccc}
\K {\mathcal A}(U\cap V) & \rightarrow & \K {\mathcal A}(U) \\
\downarrow & & \downarrow \\
\K {\mathcal A}(V) & \rightarrow & \K {\mathcal A}(X) \\
\end{array} . \]

Finally, note that the disjoint union property follows immediately from the definition of the algebroid
${\mathcal A} (X)$, and the fact that algebraic $K$-theory is compatible with direct sums and direct limits.
\end{proof}

The above result, along with theorem \ref{CH}, gives us a new proof of the fact that the functor
$X\mapsto \K ({\mathcal A}({\mathcal O}X))$ is excisive.\footnote{For other proofs
of this fact, see for instance \cite{ACFP,Weiss}.}

\begin{corollary} \label{Ahom2}
The set of functors $\{ X\mapsto K_n {\mathcal A} (X) \ |\ n\in \Z \}$
on the category of coarse spaces is a coarse homology theory
satisfying the disjoint union axiom.
\noproof
\end{corollary}

As observed in \cite{HPR,Mitch4}, one can construct a coarse
homology theory by taking the analytic $K$-theory of the Roe
$C^\ast$-algebra of a coarse space.  However, this particular
construction needs to be modified to make something that is
functorial at the level of spectra.

Let $\mathcal A$ be a $\mathbb C$-algebroid.  Then we call $\mathcal
A$ a {\em Banach category} if each morphism set $\Hom
(a,b)_{\mathcal A}$ is a Banach space, and given morphisms $x\in
\Hom (a,b)_{\mathcal A}$ and $y\in \Hom (b,c)_{\mathcal A}$, we have
the inequality $\| yx \| \leq \| y\| \cdot \| x\|$.

An {\em involution} on a Banach category $\mathcal A$ is a
collection of conjugate-linear maps $\Hom (a,b)_{\mathcal A}
\rightarrow \Hom (b,a)_{\mathcal A}$, written $x\mapsto x^\ast$,
such that $(x^\ast )^\ast =x$ for every morphism $x$, and
$(xy)^\ast = y^\ast x^\ast$ whenever $x$ and $y$ are composable
morphisms.

A Banach category with involution is termed a {\em $C^\ast$-category} if
for every morphism $x\in \Hom (a,b)_{\mathcal A}$ the {\em
$C^\ast$-identity} $\| x^\ast x \| = \| x\|^2$ holds, and the
element $x^\ast x$ is a positive element of the $C^\ast$-algebra
$\Hom (a,a)_{\mathcal A}$.

We refer the reader to \cite{GLR,Mitch2} for more information on the
theory of $C^\ast$-categories.  As in \cite{Jo2,Mitch2.5}, one can
associate a $K$-theory spectrum ${\mathbb K} ({\mathcal A})$ to a
$C^\ast$-category $\mathcal A$.  The stable homotopy groups of the
$K$-theory spectrum of a $C^\ast$-algebra are the usual analytic
$K$-theory groups of a $C^\ast$-algebra.

A functor between $C^\ast$-categories that is linear on each morphism set and compatible with the involution is called a {\em
$C^\ast$-functor}.  A $C^\ast$-functor $\alpha \colon {\mathcal
A}\rightarrow {\mathcal B}$ is automatically continuous on each
morphism set, with norm at most one, and induces a map
$\alpha_\ast \colon \K ({\mathcal A}) \rightarrow \K ({\mathcal B})$
of $K$-theory spectra.  With these induced maps, formation of the
$K$-theory spectrum defines a functor, $\K$, from the category of
$C^\ast$-categories and $C^\ast$-functors to the category of
spectra.

The functor $\K$ takes short exact sequences
of $C^\ast$-categories to fibrations of spectra, and takes
equivalences of $C^\ast$-categories to homotopy-equivalences of
spectra.

The following definition is inspired by \cite{HP}, though it is more
general.

\begin{definition}
Let $\mathcal A$ be an additive $C^\ast$-category.  Then we define
${\mathcal A}^b (X)$ to be the category of geometric $\mathcal
A$-modules over $X$, and morphisms $\phi \colon M\rightarrow N$ such
that the linear map
\[ T_\phi \colon \oplus_{x\in X} M_x \rightarrow \oplus_{x\in X}N_x \]
defined by the formula
\[ T_\phi (v) = \sum_{y\in X} \phi_{x,y} (v) \qquad v\in M_x \]
is bounded.  In this case, we define the {\em norm} of the morphism $\phi$ by writing $\| \phi \| = \| T_\phi \|$.
\end{definition}

Observe that ${\mathcal A}^b (X)$ is a {\em pre-$C^\ast$-category}
in the sense that it has all of the properties required of a
$C^\ast$-category apart from the morphism sets being complete.  As
explained in \cite{Mitch2}, we can therefore complete it to form a
$C^\ast$-category, which we label ${\mathcal A}^\ast (X)$.

The following result can now be proved in the same way as theorem
\ref{Ahom}.

\begin{theorem} \label{Ahom3}
Let $\mathcal A$ be a $C^\ast$-category.  Then
the functor $X\mapsto {\mathbb K} {\mathcal A}^\ast (X)$ is coarsely
excisive.
\noproof
\end{theorem}

The following result is proved in \cite{Mitch12}, and relates the above construction to the Roe $C^\ast$-algebra, $C^\ast (X)$.  The $K$-theory of the Roe $C^\ast$-algebra is the right-hand side of the coarse Baum-Connes conjecture, as discussed in \cite{HR2}, and can be considered the prototype of the other coarse isomorphism conjectures mentioned in this article.

\begin{theorem} \label{NEWc}
Let $\mathcal V$ be the $C^\ast$-category where the objects are
the Hilbert spaces $\C^n$, and the morphisms are bounded linear
maps.  Let $C^\ast (X)$ be the Roe $C^\ast$-algebra of the coarse topological space $X$.

Then we have a canonical
$C^\ast$-functor $\alpha \colon {\mathcal V}^\ast (X)\rightarrow
C^\ast (X)$ that induces an equivalence of
$K$-theory spectra.
\noproof
\end{theorem}

In particular, the assignment $X\mapsto {\mathbb K}C^\ast (X)$ is functorial; this fact is far from obvious if looked at directly (although the $K$-theory groups of the Roe $C^\ast$-algebra are easily seen to be functors).

\section{Equivariant Assembly}

Let $G$ be a discrete group.

\begin{definition}
A coarse space equipped with an action of the group $G$
by coarse maps is termed a {\em coarse $G$-space}.
\end{definition}

We will assume that a group $G$ acts on the right of a space.  We call a subset, $A$,
of a coarse $G$-space $X$ {\em cobounded} if there is a bounded subset $B\subseteq X$
such that $A\subseteq BG$.

\begin{definition}
The {\em coarse $G$-category} is the category where the objects are coarse $G$-spaces, and the morphisms
are controlled equivariant maps where the inverse image of a cobounded set is cobounded.
\end{definition}

If $X$ is a coarse $G$-space, and $p\colon X\rightarrow R$ is a map to a generalised
ray, the group $G$ acts on the cylinder $I_p X$ by writing $(x,t)g =
(xg,t)$.  The inclusions $i_0,i_1 \colon X\rightarrow I_p X$ are morphisms in the coarse $G$-category.

\begin{definition}
Let $f_0,f_1\colon X\rightarrow Y$ be morphisms between
coarse $G$-spaces.  An {\em elementary coarse $G$-homotopy} between
$f_0$ and $f_1$ is an equivariant coarse map
$H\colon I_p X\rightarrow Y$ for some $p\colon X\rightarrow R$ such that $f_0 = H\circ i_0$
and $f_1 = H\circ i_1$.

More generally, we call the maps $f_0$ and $f_1$ {\em coarsely $G$-homotopic} if they can be linked by a chain of elementary coarse
homotopies.
\end{definition}

A morphism $f\colon X\rightarrow Y$ in the coarse $G$-category is termed a {\em coarse
$G$-homotopy equivalence} if there is a morphism $g\colon Y\rightarrow
X$ such that the compositions $g\circ f$ and $f\circ g$ are coarsely
$G$-homotopic to the identities $1_X$ and $1_Y$ respectively.

\begin{definition}
We call a coarse $G$ space $X$ {\em $G$-flasque} if there is an equivariant map
$t\colon X\rightarrow X$ such that:

\begin{itemize}

\item Let $B\subseteq X$ be bounded.  Then there exists $N\in \N$
such that $t^n[X]\cap B = \emptyset$ for all $n\geq N$.

\item Let $M\subseteq X\times X$ be controlled.  Then the union
$\bigcup_{n\in \N} t^n [M]$ is controlled.

\item The map $t$ is close to the identity map.

\end{itemize}

\end{definition}

Note that the above map $t$ is controlled by the above axioms.

As in the non-equivariant case, if a coarse $G$-space $X$ is $G$-flasque, then
so is the closed cone ${\mathcal C}X$.  If $X$ is any coarse
$G$-space, the open square ${\mathcal S}X$ is $G$-flasque.

\begin{definition}
We call a functor $E_G$ from the coarse $G$-category to the category
of spectra {\em coarsely $G$-excisive} if the following conditions
hold.

\begin{itemize}

\item The spectrum $E_G(X)$ is weakly contractible whenever the coarse
space $X$ is $G$-flasque.

\item The functor $E_G$ takes coarse $G$-homotopy equivalences to weak homotopy
equivalences of spectra.

\item Given a coarsely excisive decomposition $X=A\cup B$, where $A$ and $B$
are coarse $G$-spaces, we have a homotopy push-out diagram
\[ \begin{array}{ccc}
E_G(A\cap B) & \rightarrow & E_G(A) \\
\downarrow & & \downarrow \\
E_G(B) & \rightarrow & E_G(X) \\
\end{array} . \]

\item Up to homotopy, the functor $E_G$ takes disjoint unions in the coarse $G$-category to coproducts
in the category of spectra.

\item Let $X$ be a cobounded coarse $G$-space.  Then the constant map $c\colon X\rightarrow +$ induces a stable
equivalence $c_\ast \colon E_G(X)\rightarrow E_G(+)$.

\end{itemize}

\end{definition}

Our next definition comes from \cite{DL}.

\begin{definition}
Let $\mathcal C$ be the category of topological spaces that are
homotopy-equivalent to $CW$-complexes.  We call a functor, $F$, from
the category $\mathcal C$ to the category of spectra {\em $G$-excisive} if it takes $G$-homotopy equivalences of spaces to homotopy
equivalences of spectra, up to homotopy takes disjoint unions of $G$-spaces
to coproducts in the category of spectra, and given a space
$X=U\cup V$, where $U$ and $V$ are $G$-invariant open sets we have a homotopy
push-out diagram
\[ \begin{array}{ccc}
F(U\cap V) & \rightarrow & F(U) \\
\downarrow & & \downarrow \\
F(V) & \rightarrow & F(X) \\
\end{array} . \]
\end{definition}

Recall that we call a $G$-space $X$ {\em cocompact} if there is a compact subset $K\subseteq X$ such
that $X=KG$.  The following is proved similarly to theorem \ref{CH}

\begin{theorem} \label{Gexc}
Let $E_G$ be a coarsely $G$-excisive functor.  Then the assignment $X\mapsto E_G({\mathcal O}X)$ defines a
$G$-excisive functor on the category of cocompact Hausdorff $G$-spaces.
\noproof
\end{theorem}

In fact, we could go further and show that the collection of functors $X\mapsto \pi_n E_G ({\mathcal O}X)$ for different groups $G$ defines an equivariant homology theory in the sense of \cite{KL}, but we do not need this stronger result here.

The following is similar to lemma \ref{csefib}.

\begin{lemma} \label{csefibg}
Let $X$ be a coarse Hausdorff $G$-space.  Let $E_G$ be a coarsely $G$-excisive
functor.
Then we have a weak fibration
\[ E_G(X) \stackrel{j}{\rightarrow} E_G({\mathcal C}X) \stackrel{v_G}{\rightarrow} E_G({\mathcal O}X) \]
\noproof
\end{lemma}

So we have an associated boundary map
\[ \partial_G \colon \Omega E_G({\mathcal O}X ) \rightarrow E_G(X) .  \]

We call this map the {\em equivariant assembly map} associated to the functor $E_G$.

Recall that a $G$-space, $X$, is termed {\em free} if for every point $x\in X$, we have $xg=x$ only when $g=e$, where $e$ denotes 
the identity element of a group $G$.  We define $EG$ to be a weakly contractible free $G$-$CW$-complex.  The space $EG$ is
unique up to $G$-homotopy-equivalence, and the quotient space $EG/G$ is the classifying space $BG$.

\begin{definition}
The {\em Novikov conjecture} associated to the functor
$E_G$ asserts that the equivariant assembly map is injective at the level of stable homotopy groups when $X=EG$ for some coarse structure on $EG$ compatible with the topology.
\end{definition}

More generally, let $\mathcal F$ be a family of subgroups of a group $G$ (that is a collection of subgroups closed under the 
operations of conjugation and finite intersections).  For instance, we can consider all finite subgroups, or all virtually cyclic subgroups.  We call a $G$-$CW$-complex $X$ a {\em $(G,{\mathcal F})$-$CW$-complex} if the fixed point set $X^H$ is empty if $H\not\in {\mathcal F}$.

Each cell in a $(G,{\mathcal F})$-$CW$-complex is a $G$-space, $C$, with the property that $C^H=\empty$ if $H\not\in {\mathcal F}$.

We can form a unique (up to $G$-homotopy) $(G,{\mathcal F})$-$CW$-complex $E(G,{\mathcal F})$, with the property that the fixed point set
$E(G,{\mathcal F})^H$ is $G$-contractible if $H\in {\mathcal F}$.

As a special case, if the family $\mathcal F$ consists of just the trivial subgroup of $G$, then $E(G,{\mathcal F}) = EG$.  If the family $\mathcal F$ consists of all finite subgroups, then $E(G,{\mathcal F})$ is the classifying space $\underline{E}G$ described for instance in \cite{BCH}.

It is shown in \cite{DL} that the space $E(G,{\mathcal F})^H$ is a classifying space for $(G,{\mathcal F})$-$CW$-complexes in the following sense.

\begin{theorem}
Let $\mathcal F$ be a family of subgroups of a group $G$.  Let $X$ be a $(G,{\mathcal F})$-$CW$-complex.  Then we have an equivariant map $u\colon X\rightarrow E(G,{\mathcal F})$, and any two such maps are $G$-homotopic.
\noproof
\end{theorem}

\begin{definition}
The $(E_G, {\mathcal F})$-isomorphism conjecture asserts that the equivariant assembly map $\partial_G \colon \Omega E_G({\mathcal O}X ) \rightarrow E_G(X)$
is a stable equivalence when $X= E(G,{\mathcal F})$ for some coarse structure on $E(G,{\mathcal F})$ compatible with the topology.
\end{definition}

We will identify examples of these conjectures in the next section.

\begin{definition}
Let $E$ be a coarsely excisive functor.  Then we say a coarsely $G$-excisive functor $E_G$ has the {\em local property} relative to $E$ if there is a natural
map $i \colon E_G(X)\rightarrow E(X)$, such that if $X={\mathcal O}Y$, where $Y$ is a free coarse cocompact $G$-space, and $\pi \colon X\rightarrow X/G$ is the quotient map, then the composite
\[ \pi_\ast \circ i = i \circ \pi_\ast \colon E_G (X)\rightarrow E(X/G) \]
is a stable equivalence.
\end{definition}

Thus, if the coarsely $G$-excisive functor $E_G$ is associated to the coarsely excisive functor $E$, we have a stable equivalence $E_G(X)\rightarrow E(X/G)$.  In particular, if
the functor $E_G$ has the local property relative to $E$, then the
assembly map in the Novikov conjecture can be regarded as a map
\[ \partial_G \colon \Omega E({\mathcal O}BG) \rightarrow E_G(EG) . \]

\begin{proposition}
Let $E_G$ be a coarsely $G$-excisive functor with the local property relative to $E$.  Suppose that the assembly map in the $(E,G,{\mathcal F})$-isomorphism conjecture is rationally injective at
the level of stable homotopy groups for some family $\mathcal F$, and the spaces $EG$ and $E(G,{\mathcal F})$ are both cocompact.  Then the Novikov conjecture holds for the functor $E_G$.
\end{proposition}

\begin{proof}
Certainly, the space $EG$ is a $(G,{\mathcal F})$-$CW$-complex, so we have a map $u\colon EG\rightarrow E(G,{\mathcal F})$.  On the other hand,
the space $E(G,{\mathcal F})$ is a $G$-space, so we have a map $v\colon E(G,{\mathcal F})\rightarrow BG$.  Since $EG$ is a $G$-space, up to $G$-homotopy there is
only one equivariant map $EG\rightarrow BG$.  But $BG =EG/G$, so both $v\circ u$ and the quotient map $\pi$ are such maps.

We conclude that the quotient map $\pi$ and the composite $v\circ u$ are $G$-homotopic.  Since the functor $E_G$ has the local property, and $EG$ is a free
$G$-space, we see that the composite $v_\ast \circ u_\ast \colon E_G ({\mathcal O}EG)\rightarrow E_G ({\mathcal O}BG)$ is a stable equivalence.  In particular, the functor $u_\ast \colon E_G ({\mathcal O}EG)\rightarrow E_G( {\mathcal O} E(G,{\mathcal F}))$ is split injective at the level of stable homotopy groups.

Now, the map $u$ induces a commutative diagram
\[ \begin{array}{ccc}
\Omega E_G({\mathcal O}EG ) & \rightarrow & E_G(EG) \\
\downarrow & & \downarrow \\
\Omega E_G({\mathcal O}E(G,{\mathcal F}) ) & \rightarrow & E_G(E(G,{\mathcal F})) \\
\end{array} . \]

By the last of the axioms involved in the definition of a coarsely excisive functor, the spaces $E_G(EG)$ and $E_G(E(G,{\mathcal F}))$ are both
naturally stably equivalent to the space $E_G(+)$, so the downward arrow on the right in the above diagram is a stable equivalence.

We saw above that the vertical arrow on the left is split injective at the level of stable homotopy groups, and by hypothesis, the map at the bottom is injective.  So the map at the top is also injective at the level of stable homotopy groups.  But this statement is the Novikov conjecture.
\end{proof}

\section{Equivariant Examples}

Let $X$ be a coarse $G$-space, let $R$ be a ring, and let $\mathcal A$ be an additive $R$-algebroid.  Then we call a geometric $\mathcal A$-module, $M$,
over $X$ {\em $G$-invariant} if $M_{xg} = M_x$ for all $x\in X$ and $g\in G$. A morphism $\phi \colon M\rightarrow N$ between such modules is termed {\em $G$-invariant}
if $\phi_{xg,yg} = \phi_{x,y}$ for all $x,y\in X$.

\begin{definition}
We write ${\mathcal A}_G [X]$ to denote the category of $G$-invariant geometric $\mathcal A$-modules over $X$, and $G$-invariant morphisms.
\end{definition}

The following result is similar to theorem \ref{Ahom}.

\begin{theorem}
The assignment $X\mapsto \K {\mathcal A}_G [X]$ is a coarsely $G$-excisive functor.
\noproof
\end{theorem}

The following is obvious from the definition.

\begin{proposition}
The functor $\K {\mathcal A}_G$ has the local property relative to the functor ${\mathcal A}$.
\noproof
\end{proposition}

Of course, it follows from the above that the functor $X\mapsto \K {\mathcal A}_G [{\mathcal O}X]$ is $G$-excisive.  An alternative proof of this can be found for instance in \cite{BFJR}.

\begin{definition}
Let $X$ be a coarse $G$-space, and let $\mathcal A$ a $C^\ast$-category.  Then we define ${\mathcal A}_G^b (X)$
to be the category of $G$-invariant geometric $\mathcal A$-modules over $X$, and $G$-invariant bounded morphisms.

We define ${\mathcal A}_G^\ast (X)$ to be the $C^\ast$-category we obtain by completion.
\end{definition}

Similarly to theorem \ref{Ahom3}, we have the following

\begin{theorem} \label{Ahom3}
Let $\mathcal A$ be a $C^\ast$-category.  Then
the functor $X\mapsto {\mathbb K} {\mathcal A}_G^\ast (X)$ is coarsely $G$-excisive.
\noproof
\end{theorem}

Again, the functor $\K {\mathcal A}_G^\ast$ has the local property relative to the functor $\K {\mathcal A}^\ast$.

Now (see for instance \cite{HR1}), given a coarse topological space $X$, there is an equivariant analogue of the Roe $C^\ast$-algebra, $C^\ast_G (X)$.

Let $|G|$ be the coarse space associated to a finitely presented group $G$, by picking a word length metric.\footnote{The particular choice of word length metric does not affect the coarse structure.}  Then it is shown in \cite{Roe5} that the $C^\ast$-algebra $C^\ast_G |G|$ has the same $K$-theory as the reduced group $C^\ast$-algebra $C^\ast_r G$.

Just as in the non-equivariant case, we have the following.

\begin{theorem}
We have a $C^\ast$-functor $\alpha \colon {\mathcal V}^\ast_G (X)\rightarrow
C^\ast_G (X)$ that induces an isomorphism between
$K$-theory groups.
\noproof
\end{theorem}

Now, fix an algebroid $\mathcal A$, and consider the $G$-excisive functor defined by the formula
\[ F_G (X) = \Omega \K {\mathcal A}_G {\mathcal O}(X) . \]

According to \cite{DL}, the Farrell-Jones assembly map in algebraic $K$-theory is the map $c\colon F_G (X)\rightarrow F_G( +)$
induced by the constant map $X\rightarrow +$.  The following result closely follows the methods described in \cite{HP} for controlled assembly.

\begin{theorem}
Let $G$ act cocompactly on the space $X$.  Then the Farrell-Jones assembly map $c\colon F_G(X)\rightarrow F_G(+)$ is stably-equivalent to the
equivariant assembly map $\partial_G \colon \Omega \K {\mathcal A}_G ({\mathcal O}X) \rightarrow \K {\mathcal A}_G(X)$.
\end{theorem}

\begin{proof}
Consider the commutative diagram
\[ \begin{array}{ccccc}
\K {\mathcal A}_G (X) & \rightarrow & \K {\mathcal A}_G ({\mathcal C}X ) & \rightarrow & \K {\mathcal A}_G ({\mathcal O}X) \\
\downarrow & & \downarrow & & \downarrow \\
\K {\mathcal A}_G (+) & \rightarrow & \K {\mathcal A}_G ({\mathcal C}+ ) & \rightarrow & \K {\mathcal A}_G ({\mathcal O}+) \\
\end{array} \]
where the vertical maps are all induced by the constant map $X\rightarrow +$.  Since the $G$-space $X$ is cocompact, it is cobounded
as a coarse $G$-space, so constant map $X\rightarrow +$ is a morphism in the coarse $G$-category.

The rows of the above diagram are fibrations.  The space ${\mathcal C}+$ is clearly flasque, so the spectrum $\K {\mathcal A}_G ({\mathcal C}+)$ is weakly contractible.  We therefore have a commutative
diagram
\[ \begin{array}{ccc}
\Omega \K {\mathcal A}_G ({\mathcal O}X ) & \rightarrow & \K {\mathcal A}_G (X) \\
\downarrow & \downarrow & \downarrow \\
\Omega \K {\mathcal A}_G ({\mathcal O}+ ) & \rightarrow & \K {\mathcal A}_G (+) \\
\end{array} \]
where the bottom row is a stable equivalence.

Since the group $G$ acts cocompactly on the space $X$, there is a bounded set $K$ such that $GK =X$, and $xg\not\in K$ if $x\in K$ and $g\neq e$.  Thus, if we give the space
$K$ the trivial $G$-action, by definition of the functor ${\mathcal A}_G$, the categories ${\mathcal A}_G (X)$ and ${\mathcal A}_G (K)$ are equivalent.

But the spaces $K$ and $+$ are equivariantly coarsely equivalent, and so equivariant coarse homotopy-equivalent.  Hence the map $a \colon \K {\mathcal A}_G [X]\rightarrow \K {\mathcal A}_G [+]$
is a stable equivalence.

The desired result now follows.
\end{proof}

The corresponding result for the Baum-Connes assembly map follows similarly, although we need the fact that the assembly map fits into the picture in \cite{DL}; see for example \cite{Mitch6}.

\begin{theorem}
Let $G$ act cocompactly on the space $X$.  Then the Baum-Connes assembly map for the $G$-space $X$ is stably-equivalent to the
equivariant assembly map $\partial_G \colon \Omega \K {\mathcal V}^\ast_G ({\mathcal O}X) \rightarrow \K {\mathcal V}^\ast_G(X)$.
\noproof
\end{theorem}

\section{Descent}

Let $G$ be a discrete group, and let $A$ be a $G$-spectrum.   Then we
define the {\em homotopy fixed point spectrum}, $A^{hG}$, to be the
spectrum of continuous equivariant maps $\mu \colon EG\rightarrow
Y$.  Each space in the spectrum $A^{hG}$ has the compact open topology.

The following result is well-known.

\begin{proposition} \label{WK}
Let $A$ and $B$ be $G$-spectra.  Let $f\colon A\rightarrow B$ be a (non-equivariant) weak homotopy-equivalence.  Then the obvious induced map $f_\ast \colon A^{hG}\rightarrow B^{hG}$ is also a weak homotopy-equivalence.

In particular, if the spectrum $A$ is weakly contractible, then so is the spectrum $A^{hG}$.
\end{proposition}

A proof of the second assertion (at least for spaces) can be found in \cite{Roe1}; the first then follows by looking at mapping cones.

Applying the above to a coarsely excisive functor immediately gives us the following.

\begin{proposition} \label{EHG}
Let $E$ be a coarsely excisive functor.  Then the functor $X\mapsto
E(X)^{hG}$, defined on the category of coarse $G$-spaces and equivariant
coarse maps, has the following three properties.
 
\begin{itemize}

\item Let $X$ be $G$-flasque.  Then the spectrum $E(X)^{hG}$ is weakly contractible.

\item Let $f\colon X\rightarrow Y$ be an equivaraint coarse homotopy equivalence.  Then the induced map $f_\ast \colon E(X)^{hG}\rightarrow
E(Y)^{hG}$ is a weak homotopy equivalence.

\item Let $X=A\cup B$ be a coarsely excisive decomposition $X=A\cup B$, where $A$ and $B$
are coarse $G$-spaces.  Then we have a homotopy push-out
diagram
\[ \begin{array}{ccc}
E(A\cap B)^{hG} & \rightarrow & E(A)^{hG} \\
\downarrow & & \downarrow \\
E(B)^{hG} & \rightarrow & E(X)^{hG} \\
\end{array} . \]

\end{itemize}

\noproof
\end{proposition}

Our next result can be proved using the above in a similar way to theorem \ref{CH}

\begin{corollary}
The sequence of functors $X\mapsto \pi_n E({\mathcal O}X)^{hG}$ forms a $G$-homology theory.
\noproof
\end{corollary}

By a $G$-homology theory, we mean one that satisfies equivariant analogues of the Eilenberg-Steenrod axioms; see \cite{KL} for
a discussion.  Note, however, that the disjoint union axiom need not hold in this case.

Proposition \ref{EHG} is strong enough for the following result to hold.  The proof is similar to lemma \ref{csefib}.

\begin{corollary}
Let $E$ ba a coarsely excisive functor, and let $X$ be a coarse
Hausdorff $G$-space.  Then we have a weak fibration
\[ E(X)^{hG} \rightarrow E({\mathcal C}X)^{hG}
\stackrel{v_{hG}}{\rightarrow} E({\mathcal O}X)^{hG} \]
\noproof
\end{corollary}

So we have a boundary map
\[ \partial_{hG} \colon \Omega E({\mathcal O}X)^{hG} \rightarrow
E(X)^{hG} \]

\begin{lemma}
Let $E_G$ be a coarsely $G$-excisive functor with the local property relative to a functor $E$.  Then we have a natural transformation $j\colon E_G (X)\rightarrow E(X)^{hG}$.

Further, let $X=\vee_{g\in G}Y_g$, where each space $Y_g$ is a copy of the same coarse space $Y$, which is a cone of some compact space, and the group $G$
acts by permutations.  Then the map $j\colon E_G (X)\rightarrow E(X)^{hG}$ is an isomorphism.
\end{lemma}

\begin{proof}
By hypothesis, we have a natural map $i\colon E_G (X)\rightarrow E(X)$.

Observe we have a natural homotopy equivalence
\[ \begin{array}{rcl} 
E_G (X) & \simeq & \Map (EG,E_G(X)) \\
& = & \Map (EG,\Map_G(G,E_G(X))) \\
& = & \Map_G (EG,\Map (G,E_G(X))) \\
& = & \Map (G,E_G(X))^{hG} 
\end{array}. \]

We have a natural map $\Map (G,E_G(X)) \rightarrow E(X)$ defined by composing the evaluation of a map $f\colon G\rightarrow E_G(X)$ at the identity element $1\in G$, with the natural map $j$.

Taking homotopy fixed point sets, we have a natural transformation $i\colon E_G(X)\rightarrow E(X)^{hG}$.

Now, let $X=\vee_{g\in G}Y_g$ as above.  Let $Y= {\mathcal O}Z$.  Then $X= {\mathcal O}(Z\times G)$.   

Let $\pi \colon X\rightarrow Y$ be the quotient map.  Then by the local property, the map $\pi_\ast \circ i = i \circ \pi_\ast \colon E_G(X)\rightarrow E(Y)$ is a stable equivalence.  Now, looking at disjoint unions,
\[ E(X) = \vee_{g\in G} E(Y_g) = \Map (G, E(Y) ) . \]

So in this case tha above map $\Map (G,E_G(X))\rightarrow \Map (G,E(Y))$ is weakly homotopic to the stable equivalence $\Map (G,E_G(X))\rightarrow \Map (G,E(Y))$ induced by the composite $\pi_\ast \circ i$.  In particular, it follows that our map is a stable equivalences.

Taking homotopy fixed points, the map
\[ j\colon E_G(X) \rightarrow E(X)^{hG} \]
is also a stable equivalence, and we are done.
\end{proof}

So we have natural maps $i_1 \colon E_G({\mathcal C}X)\rightarrow E({\mathcal C}X)^{hG}$ and
$i_2\colon E_G({\mathcal O}X) \rightarrow E({\mathcal O}X)^{hG}$
fitting into a commutative diagram
\[ \begin{array}{ccc}
E_G({\mathcal C}X) & \stackrel{v_G}{\rightarrow} & E_G({\mathcal
O}X) \\
\downarrow & & \downarrow \\
E({\mathcal C}X)^{hG} & \stackrel{v_{hG}}{\rightarrow} & E({\mathcal
O}X)^{hG} \\
\end{array} \]

\begin{proposition} \label{ihg}
The map $i_2 \colon E_G ({\mathcal O}X)\rightarrow E({\mathcal O}X)^{hG}$ is a stable
equivalence whenever $X$ is a finite cocompact free $G$-$CW$-complex.
\end{proposition}

\begin{proof}
Let $h_\ast^G$ be a $G$-homology theory.  Let $X$ be a finite free
$G$-$CW$-complex.  Then the equivariant version
of the Atiyah-Hirzebruch spectral sequence (see for instance \cite{Mat}) gives us a half-plane spectral sequence, $\{ E^\ast_{p,q} \}$,
converging to $h_\ast^G (X)$, with $E^2$-term $E^2_{p,q} \cong H_p^G (X;h_q^G (\Delta^0 ) )$, that is classical $G$-homology with coefficients in the group $h_q^G (\Delta_G^0 )$,
where $\Delta_G^0$ is the free $0$-dimensional $G$-cell.  

Further, this spectral sequence depends functorially on the $CW$-complex $X$ and $G$-homology theory $h_\ast^G$.  The space $\Delta_G^0$ is just a copy of the group $G$, which acts on itself by right-translations.

Hence, the induced map $(i_2)_\ast \colon \pi_\ast E_G({\mathcal O}X) \rightarrow E({\mathcal O}X)^{hG}$ is a map between $G$-homology theories.  By the above, it is an isomorphism
when the space $X$ is a $0$-dimensional free $G$-cell.  Therefore the corresponding Atiyah-Hirzebruch spectral sequences are isomorphic, and the map $(i_2)_\ast$ is an isomorphism for any finite free $G$-$CW$-complex.
\end{proof}

We now have the ingredients to prove our main result.

\begin{theorem}
Let $E_G$ be a coarsely $G$-excisive functor.  Let $X$ be a free coarse $G$-space, that is, as a topological space, $G$-homotopy equivalent to a finite $G$-$CW$-complex.

Suppose the coarse isomorphism conjecture holds for the functor $E$ and the space $X$.  Then the map $\partial_G \colon E_G ({\mathcal O}X) \rightarrow E_G (X)$ is injective at the level of stable homotopy groups.
\end{theorem}

\begin{proof}
By proposition \ref{wcx}, the coarse isomorphism conjecture for the space $X$ tells us
that the spectrum $E({\mathcal C}X)$ is weakly contractible.

Hence, by proposition \ref{WK}, the spectrum $E({\mathcal C}X)^{hG}$ is also weakly contractible.

Now we have a commutative diagram
\[ \begin{array}{ccc}
E_G({\mathcal C}X) & \stackrel{v_G}{\rightarrow} & E_G({\mathcal
O}X) \\
\downarrow & & \downarrow \\
E({\mathcal C}X)^{hG} & \stackrel{v_{hG}}{\rightarrow} & E({\mathcal
O}X)^{hG} \\
\end{array} \]
and by proposition \ref{ihg} the map $i_2 \colon E_G ({\mathcal O}X)\rightarrow E({\mathcal O}X)^{hG}$
is a weak equivalence.

Hence the map $v_G$ must be zero at the level of stable homotopy
groups.  The weak fibration in lemma \ref{csefibg} gives us a long
exact sequence
\[ \pi_{n+1} E({\mathcal C}X)^G \stackrel{(v_G)_\ast }{\rightarrow} \pi_{n+1} E({\mathcal O}X)^G
\stackrel{(\partial_G)_\ast }{\rightarrow}  \pi_n E(X)^G \]

As we have just mentioned, the map $(v_G)_\ast$ is zero.  Thus the
map $(\partial_G)_\ast$ is injective. 
\end{proof}

Applying the above result to the examples in the previous section, we immediately obtain the following.

\begin{corollary}
Let $X$ be a free coarse $G$-space, that is, as a topological space, $G$-homotopy equivalent to a finite $G$-$CW$-complex.

Let $R$ be a ring, and let $\mathcal A$ be an additive $R$-algebroid.  Suppose the coarse assembly map $\partial \colon \Omega \K{\mathcal A}({\mathcal O}X)\rightarrow \K {\mathcal A}(X)$ is an isomorphism.  Then the Farrell-Jones assembly map is injective for the space $X$, $R$-algebroid $\mathcal A$ and group $G$.
\noproof
\end{corollary}

\begin{corollary}
Let $X$ be a free coarse $G$-space, that is, as a topological space, $G$-homotopy equivalent to a finite $G$-$CW$-complex.

Suppose the coarse assembly map $\partial \colon \Omega \K{\mathcal V}^\ast({\mathcal O}X)\rightarrow \K {\mathcal V}^\ast (X)$ is an isomorphism.  Then the Baum-Connes assembly map is injective for the space $X$ and group $G$.
\noproof
\end{corollary}

\bibliographystyle{plain}

\end{document}